\input amstex
\documentstyle{amsppt}

\magnification=1200

\def\k{\Bbbk}
\def\Eij{E_{ij}}
\def\gln{\frak{gl}(n)}
\def\glN{\frak{gl}(N)}
\def\Ugln{\Cal U(\gln)}

\def\Zgln{\frak Z(\gln)}
\def\ZglN{\frak Z(\glN)}
\def\dt{\delta}

\def\tht{\thetag}
\def\C{\Bbb C}
\def\Cp{\bold C}
\def\kSk{\k[S(k)]}
\def\A{\Bbb A}
\def\S{\Bbb S}
\def\Sm{\S_{\mu}}
\def\nS{\overline{\S}}
\def\Z{\Bbb Z}
\def\tr{\operatorname{tr}}

\def\sgn{\operatorname{sgn}}
\def\l{\lambda}
\def\chm{\chi^\mu}
\def\s{s^*}
\def\sm{\s_{\mu}}
\def\Alt{\operatorname{Alt}}

\def\P{\Cal P}
\def\Q{\Cal Q}

\def\Tt{\widetilde{T}}
\def\la{\left\langle}
\def\ra{\right\rangle}
\def\M{\operatorname{M}}
\def\Ind{\operatorname{Ind}}
\def\Res{\operatorname{Res}}
\def\k{\Bbbk}
\def\f{\downharpoonright}
\def\u{\upharpoonright}
\def\RTab{\operatorname{RTab}}
\def\a{\alpha}
\def\Ls{\Lambda^*}
\def\Tr{T^r}
\def\gi{(g^{-1})}
\def\Dm{\Delta_\mu}
\def\ov{\,/\,}
\def\Ad{\operatorname{Ad}}
\def\Sim{\Sigma_{\mu}}

\NoBlackBoxes

\pretitle
{
\noindent \sl Submitted to Transformation Groups
\vskip 2cm
}

\topmatter
\title
Quantum immanants \\
and higher Capelli identities
\endtitle
\author
Andrei Okounkov
\endauthor
\thanks
The author is supported by the International Science
Foundation and the Russian Fundamental Research Foundation.
\endthanks
\abstract
We consider remarkable central elements of the universal
enveloping algebra $\Ugln$ which we call quantum immanants.
We express them in terms of generators $E_{ij}$ of $\Ugln$
and as differential operators on the space of matrices.
These expressions are a direct generalization of the classical
Capelli identities. They result in many nontrivial properties
of quantum immanants.  
\endabstract
\address
Institute for Problems of Information Transmission,
Bolshoj Karetny, 19, Moscow, 101447, Russia
\endaddress
\email
okounkov\@ippi.ac.msk.su
\endemail
\endtopmatter

\head 
1. Introduction
\endhead
\subhead 1.1\endsubhead
By $\Eij$ denote the standard generators of the
universal enveloping algebra $\Ugln$.
Consider the following element  of $\Ugln$
$$
\Cp=\sum_{s\in S(n)} \sgn(s)\, E_{1,s(1)} 
(E_{2,s(2)}+\dt_{2,s(2)}) \dots 
(E_{n,s(n)}+(n-1)\dt_{n,s(n)})\,. \tag 1.1
$$
Symbolically we can write
$$
\Cp=\operatorname{row-det}
\left[
\matrix
\format \l\quad &\l\quad &\l\quad &\l \\
E_{11} & E_{12} & \hdots & E_{1n} \\
E_{21} & E_{22}+1 &  & E_{2n} \\
\vdots & & \ddots & \vdots \\
E_{n1} & E_{n2} & \hdots & E_{nn}+n-1 
\endmatrix
\right] \,, \tag 1.1${}'$
$$
where the row-determinant of this matrix with
non-commutative entries is defined by \tht{1}.

Denote by $\M(n)$ the space of $n\times n$-matrices.
Denote by $\k$ the ground field. We suppose that
$\operatorname{char}\k=0$.
Consider the representation $L$ of $\Ugln$ in
the space $\k[\M(n)]$, $i,j=1,\dots,n$ given
on the generators by the following formula
$$
L(\Eij)=\sum_\alpha x_{i\alpha} \partial_{j\alpha} \,. \tag 1.2
$$
It is well known that $L$ maps $\Ugln$ isomorphically
onto the algebra of all differential operators 
with polynomial coefficients on 
the space $\M(n)$ that commute with the right action of
$GL(n)$. Introduce formal matrices $E$, $X$, $D$ with
$(i,j)$-th entry  equal to $\Eij$, $x_{ij}$, $\partial_{ij}$
respectively. Then \tht{2} can be written as
$$
L(E)=X\cdot D' \,, \tag 1.3
$$
where prime means transposition. The celebrated Capelli identity
\cite{C} asserts that
$$
L(\Cp)= \det X \cdot \det D \,. \tag 1.4
$$
Here $\det X$ and $\det D$ are ordinary determinants.
Observe that the RHS of \tht{4} visibly commutes
with the left action of $GL(n)$ so that $\Cp$ is in
fact a central element of $\Ugln$.
The Capelli identity is one of the most
imporatant results of classical invariant
theory \cite{H}. Modern approaches to this identity
were developed in \cite{HU}, \cite{KS} and by other authors.
(See, for example, references in the cited papers.)
One of these modern approaches is 
based on the notion of a quantum determinant for Yangian $Y(\gln)$
(see \cite{MNO}).
There are a $q$-analog of the Capelli identity \cite{NUW}
and its super analog \cite{N}.

In this paper we study
some remarkable generalizations of the Capelli element 
which we call quantum immanants. In some sence 
we replace the determinant in \tht{1${}'$} and \tht{4} by the
trace of a arbitrary polynomial representation of $GL(n)$.
Our approach is based on $R$-matrix formalism (however
we do not consider Yangians).
\footnote"*"{
Recently the author proved a more general Capelli-type
identity which involves not only the center but the
whole algebra $\Ugln$. The proof (which does not require
$R$-matrices) will be given in the next paper.}

Normally the Capelli element \tht{1} is defined by 
following column determinant
$$
\Cp=\sum_{s\in S(n)} \sgn(s)\, (E_{s(1),1}+(n-1)\dt_{s(1),1}) 
(E_{s(2),2}+(n-2)\dt_{s(2),2}) \dots 
E_{s(n),n}\,. 
$$
Quantum immanants (see below) can be also rewritten
in the column form.

\subhead 1.2\endsubhead
Introduce the formal matrix
$$
E(u)=\big[\Eij- u\cdot \dt_{ij}\big]_{i,j=1}^n \,.
$$
Here $u$ is a formal variable. A formal $n\times n$ matrix $A$
with entries $a_{ij}$ from a noncommutative algebra $\A$ can
be considered as an element
$$
A=\sum_{ij} a_{ij} \otimes e_{ij}  \,\in\, \A\otimes \M(n) \,,
$$
where $e_{ij}$ are standard matrix units in $\M(n)$. The
tensor product of two such matrices $A$ and $B$ is defined by
$$
A\otimes B = \sum_{i,j,k,l} a_{ij} b_{kl}
\otimes e_{ij} \otimes e_{kl} \, \in  \A\otimes \M(n)^{\otimes 2} \,.
$$
In the space $(\k^n)^{\otimes k}$ acts the symmetric group $S(k)$
so that we have a representation
$$
\kSk\to\M(n)^{\otimes k}\,.
$$
It can be shown (in fact this is a way
to prove \tht{4}; see \cite{MNO} and below) that 
$$
\Cp=(n!)^{-1} \tr 
\big(E\otimes E(-1)\otimes \dots \otimes E(-n+1) \cdot \Alt \big) \,,
\tag 1.5
$$
where $\Alt$ is the anti-symmetrizer
$$
\Alt = \sum_{s\in S(n)} \sgn(s)\cdot s \, \in \k[S(n)] \,,
$$
and the trace of an element of $\A\otimes \M(n)^{\otimes n}$
is defined by
$$
\tr \big( \sum_
{i_1,j_1,\dots,i_n,j_n}
a_{i_1,j_1,\dots,i_n,j_n} \otimes
e_{i_1,j_1}\otimes\dots\otimes e_{i_n,j_n} \big) =
\sum_{i_1\dots,i_n}
a_{i_1,i_1,\dots,i_n,i_n} \, \in \A \,.
$$
In \tht{5} the algebra $\A$ is $\Ugln$.
In customary notations \tht{5} can be rewritten as
$$
\Cp=(n!)^{-1} \sum_{i_1,\dots,i_n} \sum_{s\in S(n)}
\sgn(s) \, E_{i_1,i_{s(1)}} \dots
(E_{i_n,i_{s(n)}}+(n-1)\dt_{i_n,i_{s(n)}}) \,.
$$
It is easy to see that the RHS of \tht{4} can be written 
in a similar form
$$
\det X \cdot \det D = 
(n!)^{-1} \tr \big(X^{\otimes n} \cdot (D')^{\otimes n}
\cdot \Alt \big) \,.
$$
Since the representation $L$ is faithful let us omit
the letter $L$ and identify elements of $\Ugln$ with
differential operators. Then the Capelli identity can
be restated as follows:
$$
\tr \big(E\otimes E(-1)\otimes \dots \otimes
E(-n+1) \cdot \Alt \big) = 
\tr \big(X^{\otimes n} \cdot (D')^{\otimes n}
\cdot \Alt \big) \,. \tag 1.6
$$
The identity \tht{6} is true also for the action
of $GL(n)$ on rectangular $n\times m$ matrices.
In this case the matrices $X$ and $D$ are also
rectangular $n\times m$ matrices. From now on
we consider this general rectangular case.

\subhead 1.3\endsubhead
Now we can formulate higher Capelli identities.
Let $\mu$ be a partition such that $\ell(\mu)\le n$.
Put $k=|\mu|$. Let $\chm$ be the character of the group $S(k)$
corresponding to the partition $\mu$
$$
\chm = \sum_{s\in S(k)} \chm(s) \cdot s \, \in \kSk \,.
$$
Let $T$ be a standard tableau of shape $\mu$. Let
$\xi_T$ be the corresponding vector of the Young
orthonormal basis. Consider the matrix element
$$
\psi_T=\sum_{s\in S(k)} (s\cdot\xi_T,\xi_T)\cdot s
\, \in \kSk \,.
$$
The element
$$
P_T=\frac{\dim \mu}{k!}\, \psi_T \, \in \kSk 
$$
acts in as
orthogonal projection onto  $\xi_T$.
in the irreducible $S(k)$-module corresponding to $\mu$
and as zero operator in other irreducible $S(k)$-modules.

Suppose $\alpha=(i,j)$ is a cell of $\mu$. The number
$c(\a)=j-i$ is called the {\it content} of the cell $\alpha$.
Suppose $\alpha$ is the $l$-th cell in the tableau $T$.
Put
$$
c_T(l)=j-i \,.
$$
For example, if
$$
T=\matrix
1 & 3 \\
2 \endmatrix 
$$
then $c_T(1)=0$, $c_T(2)=-1$, $c_T(3)=1$. 
Observe that always $c_T(1)=0$. We have

\proclaim{\smc Theorem (\sl Higher Capelli identities)} 
For all partitions $\mu$, $\ell(\mu)\le n$ and any standard tableau
$T$ of shape $\mu$
$$
\tr \big(E\otimes E(c_T(2))\otimes \dots \otimes
E(c_T(k)) \cdot P_T \big) = 
\tr \big(X^{\otimes k} \cdot (D')^{\otimes k}
\cdot \chm/k! \big) \,. \tag 1.7 
$$
In particular the LHS of \tht{7} does not depend on the 
choice of $T$.
\endproclaim

If $\mu=(1^n)$ then \tht{7} turns into \tht{6}. Below in \tht{3.24$'$}
we shall obtain another version of the identity \tht{7}
which turns into the original definition of the 
Capelli identity if $\mu=(1^n)$. 
The particular cases $\mu=(1^k)$, $k=1,\dots,n$ of this
theorem are also known as Capelli identities \cite{HU}.
A different approach to the analogs of Capelli identities for
$\mu=(k)$ can be found in \cite{N}.

\example{\smc Example} Suppose $n=1$ and $\mu=(k)$. Then
\tht{7} reads as follows
$$
x\frac{d}{dx}
\big(x\frac{d}{dx} - 1\big) \dots
\big(x\frac{d}{dx} - k+1\big) =
x^k\frac{d^k}{dx^k} \,.
$$
This identity can be easily verified by induction.
\endexample
\subhead 1.4\endsubhead
Consider the highest terms of the both sides of \tht{7}
with respect to the natural filtration in $\Ugln$.
A simple calculation shows that 
$$
\tr (E^{\otimes k} \cdot s) =
\tr (E^{\otimes k} \cdot tst^{-1}) + \text{ lower terms}
$$
for all $s,t\in S(k)$. Next observe that
$$
\sum_{t\in S(k)} t\, P_T \,t^{-1} = \chm \,.
$$
Therefore the LHS of \tht{7} equals
$$
\tr (E^{\otimes k} \cdot \chm/k!) + \text{ lower terms}.
$$
Since $X$ and $D$ commute modulo lower terms the highest
terms of the LHS and the RHS of \tht{7} agree by virtue of \tht{3}.
The structure of this highest term is similar to
the definition of the Schur function via characteristic map
\cite{M1}. Suppose $g\in GL(n)$. It follows from the classical
decomposition of $(\k^n)^{\otimes k}$ as a $GL(n)\times S(k)$-module
that the function
$$
\tr \big( g^{\otimes k} \cdot \chm/k! \big) \tag 1.8
$$
is equal to the trace of $g$ in the irreducible $GL(n)$-module
with highest weight $\mu$ or, in other words, to the Schur
polynomial in eigenvalues of $g$. Denote the polynomial \tht{8}
in matrix elements $x_{ij}$ of by $s_\mu(X)$
$$
s_\mu(X)=
\tr \big( X^{\otimes k} \cdot \chm/k! \big)\,. \tag 1.9
$$
Given a matrix $A=[a_{ij}]$, $i,j=1,\dots,k$,
the number
$$
\sum_{s\in S(k)} \chm(s) \, a_{1,s(1)} \dots a_{k,s(k)}
$$
is called the {\it $\mu$-immanant} of the matrix A.
If $\mu=(1^k),(k)$ then the $\mu$-immanant turns into
determinant and permanent respectively.
Note that \tht{9} is
the sum of $\mu$-immanants of principal
$k$-submatrices (with repeated rows and columns) of the matrix $g$.

\subhead 1.5\endsubhead
I wish to thank V.~Ginzburg, S.~Kerov, S.~Khoroshkin and M.~Noumi
for helpful discussions.
I am especially grateful to M.~Nazarov;
this paper would be hardly possible without 
numerios discussion with him.  They helped
me very much with the proof of \tht{2.2} 
(see also paragraph 2.3 below).
 
Quantum immanats and higher Capelli identities
arose from our joint work with
Olshanski \cite{OO}. The discussions we had
with G.~Olshanski during the work on \cite{OO}
were very useful for me. His critical comments
concerning this text were also very useful.

I have to mention that the structure LHS of \tht{7} is very
close to the fusion process from \cite{KuS}, \cite{KuSR},
\cite{KuR} and \cite{Ch}.

\head
2. Quantum immanants and $s^*$-functions.
\endhead

\subhead 2.1 \endsubhead
Denote the 
LHS of \tht{1.7} by $\S_\mu$:
$$
\S_\mu=
\tr \big(E\otimes E(c_T(2))\otimes \dots \otimes
E(c_T(k)) \cdot P_T \big) \,\in\Ugln \,. \tag 2.1
$$
Below we shall see that
this definition does not depend on the choice of $T$.
Because of the structure of the highest term of \tht{1}
and by analogy to the quantum determinant let us call this element
the {\it quantum $\mu$-immanant}. We shall see that quantum 
immanants have many remarkable properties.

\subhead 2.2 \endsubhead
In the next section 
we prove that the quantum immanant $\S_\mu$ lies in the
center $\Zgln$ of the algebra $\Ugln$. Next we calculate
the eigenvalue of $\pi_\l(\S_\mu)$ where $\pi_\l$
is the irreducible representation of $\Ugln$
with highest weight $\l$. We shall prove that 
$$
\pi_\l(\S_\mu)=\sm(\l) \,, \tag 2.2
$$
where $\sm$ is the {\it shifted Schur polynomial} (see \cite{OO}).
The short name of them is {\it $s^*$-polynomials}.

The definition of $s^*$-polynomials is the following.
Put
$$
\align
(a\f b)&= a(a-1)\dots(a-b+1) \\
(a\u b)&= a(a+1)\dots(a+b-1) \,.
\endalign
$$
These products are called falling and raising {\it factorial
powers}. Put also 
$$
\dt=(n-1,\dots,1,0)\,
$$
By definition \cite{OO}
$$
\sm(x_1,\dots,x_n) =\frac
{\det \big[(x_i+\dt_i\f \mu_j+\dt_j)\big]}
{\det \big[(x_i+\dt_i\f \dt_j)\big]} \,. \tag 2.3
$$
Observe that the denominator in \tht{3} equals the
Vandermode determinant in variables $x_i+\dt_i$.
Observe also that the numerator in \tht{3} is a
skew-symmetric function in $x_i+\dt_i$ and
hence is divisible by the denominator.

These polynomials were proposed by G.~Olshanski in \cite{Ol2}.
They differ by shift of variables  from the
{\it factorial\/} Schur polynomials which were
studied by L.~C.~Biedenharn and J.~D.~Louck \cite{BL},
W.~Y.~C.~Chen and J.~D.~Louck \cite{CL},
I.~Goulden and A.~Hamel \cite{GH},
I.~Goulden and C.~Greene \cite{GG},
and I.~G.~Macdonald \cite{M2}. The shift of variables is
essential. For example, in contrast to factorial
Schur polynomials the $s^*$-polynomials are stable
in the following sence
$$
s^*_\mu(x_1,\dots,x_n,0)=
s^*_\mu(x_1,\dots,x_n)\,. \tag 2.4
$$
This stability allows to introduce {\it $s^*$-functions}
in countable many variables variables as in \cite{M1}.
For example the number $\sm(\l)$, where $\l$ is a
partition, is well defined.

The $s^*$-functions have a lot of interesting
properties. Their
detailed exposition can be found in \cite{OO}.
Some of them have a natural
interpretation in terms of quantum immanants. 

\subhead 2.3 \endsubhead
Denote the differential
operator in the RHS of \tht{1.7} by $\Dm$.
We prove in the next section that $\Dm\in\Zgln$.
(Recall that we identify $\Ugln$ with the algebra of 
the right invariant differential operators.)
The higher Capelli identities \tht{1.7} will be proved in two
steps: first we prove \tht{2} and then
$$
\pi_\l(\Dm)=\sm(\l)\,. \tag 2.6
$$
The proof of \tht{6} is much more
simple than the proof of \tht{2}.

As mentioned in the introduction 
the discussions with M.~Nazarov were very
helpful for me during the proof of \tht{2}.
In particular, M.~Nazarov drew my attenton to
the importance of \tht{3.8}. He also conjectured that
the eigenvalue in $\pi_\l$ of the 
central element \tht{3.43} is equal to $\sm(\l)$.
Recently he has found a new proof of \tht{2}.

\subhead 2.4 \endsubhead
The proof of \tht{6} will be based on the two
following properties of the $s^*$-functions that
are very simple and very useful at the same time.  
The two theorems we prove below will be also
used in the forthcoming paper by A.~Molev and
M.~Nazarov concerning Capelli-type identities
for other classical groups. Many other
applications of them can be found in \cite{OO}.

The following vanishing and characterization theorems is a way
to control lower terms of inhomogeneous polynomials
$\sm$. In particular cases similar argument was
used by many people (see, for example, \cite{HU}).
In the context of Capelli identities
it was  developed in full generality by S.~Sahi in \cite{S}.
In this paper S.~Sahi considered polynomials that satisfy 
\thetag{2.7--8} and have a more general symmetry than the
shifted symmetry. He found an inductive formula for them.
However, this formula is very complicated. In our situation
we have much more simple formulas.

Denote by $H(\mu)$ the product of the hook lenghts of
all cells of $\mu$
$$
H(\mu)=\prod_{\alpha\in\mu} h(\alpha) \,.
$$
Let $\l$ be another partition $\l_1\ge\l_2\ge\dots$. Write
$\mu\subset\l$ if $\mu_i\le\l_i$ for all $i$. We have
\proclaim{\smc Vanishing theorem}
$$
\align
\sm(\l)&=0\quad\text{unless}\quad\mu\subset\lambda\,,\tag 2.7\\
\sm(\mu)&=H(\mu)\,.\tag 2.8
\endalign
$$
\endproclaim
\demo{Proof}
Observe that
$$
(a\f b)=0\quad\text{if}\quad a,b\in\Z, b>a\ge0\,.
$$
Suppose $\l_l<\mu_l$ for some $l$. Then in the matrix
$$
\big[(\l_i+n-i\f \mu_j+n-j)\big]
$$
all entries with $i\ge l$ and $j\le l$ vanish. Hence
its determinant vanishes. Since the denominator in \tht{3}
does not vanish \tht{7} follows.

Next in the matrix
$$
\big[(\mu_i+n-i\f \mu_j+n-j)\big]
$$
all entries with $i>j$ vahish. Hence its determinant equals
$$
\prod_i(\mu_i+n-i)! \,.
$$
Therefore
$$
{\sm}(\mu)=
\prod_i(\mu_i+n-i)! \,/\,
\prod_{i<j}(\mu_i-\mu_j-i+j) \,. \tag 2.9
$$
Recall that there are two formulas for the dimension
of the irreducible representation of the
symmetric group labeled by $\mu$
$$
\align
\dim\mu&=|\mu|!\,/\,H(\mu)\tag 2.10 \\
&=|\mu|!
\prod_{i<j}(\mu_i-\mu_j-i+j) \,/\,
\prod_i(\mu_i+n-i)! \,. \tag 2.11
\endalign
$$
Thus \tht{9} equals $H(\mu)$. \qed
\enddemo

\subhead 2.5 \endsubhead
By $\Ls(n)$ denote the algebra of polynomials in
variables $x_1,\dots,x_n$ which are symmetric in
new variables $x_1+\dt_1,\dots,x_n+\dt_n$. 
We call such  polynomials {\it shifted
symmetric} \cite{OO}. It is clear that $\sm\in\Ls(n)$.

Observe that the highest term of any polynomial from $\Ls(n)$
is a symmetric polynomial. It is
easy to see that the shifted Schur polynomials
$\sm$, $\ell(\mu)\le n$ form a linear basis in $\Ls(n)$.
We have
\proclaim{\smc Characterization theorem}
Any of the two following properties determines
the polinomial $\sm\in\Ls(n)$ uniquely:
\roster
\item"(A)"
$\deg {\sm}\le|\mu|$ and
$$
\sm(\l)=\dt_{\mu\l} H(\mu)
$$
for all $\l$  such that $|\l|\le|\mu|$;
\item"(B)"
the highest term of $\sm$ is the 
ordinary Schur function $s_\mu$  and
$$
\sm(\l)=0 
$$
for all $\l$ such that $|\l|<|\mu|$.
\endroster
\endproclaim

\demo{Proof} Prove part (A). 
We have to prove that
$$
\left.
\aligned
&f\in\Ls(n)\,,\\
&\deg f\le|\mu|\,,\\
&f(\l)=0\, \text{ for all }\l,\,
|\l|\le|\mu|,\, 
\ell(\l)\le n\,,
\endaligned
\right\}
\Longrightarrow
\, f=0 \,.
$$
Put $k=|\mu|$. The polynomials
$\{s^*_\l\},|\l|\le k,\ell(\l)\le n,$
is a linear basis in subspace of $\Ls(n)$ which 
consists of polynomials of degree $\le k$. Hence
$$
f=\sum c_\l s^*_\l, \quad |\l|\le k,\ell(\l)\le n\,, \tag 2.12
$$
for some coefficients $c_\l$. Show that $c_\l=0$ for
all $\l$. Suppose $c_\nu\ne0$ for some $\nu$.
Choose the partition $\nu$ so that
$c_\nu\ne0$ and $c_\eta=0$ for all $\eta$, $|\eta|<|\nu|$.
Evaluate \tht{12} at $\nu$. By the vanishing theorem we obtain
$$
0=c_\nu H(\nu) \,.
$$
Thus $c_\nu=0$.

Prove part (B).
Suppose there are two such elements $f_1$ and $f_2$ of $\Ls$.
Then $\deg(f_1-f_2)<|\mu|$ and
$(f_1-f_2)(\l)=0$ for all $\l$ such that $|\l|<|\mu|$.
By part (A) we have $f_1-f_2=0$. \qed
\enddemo

\subhead 2.6 \endsubhead
By \tht{2} the vanishing and characterization theorems can be
restated in terms of quantum immanants.
Sometimes it is convenient to use the following
normalized quantum immanants $\nS_\mu$.  Put
$$
(a\u\mu) = \prod_{\a\in\mu} (a+c(\a))\,, \tag 2.14
$$
where $c(\a)$ is the content of the cell $\a\in\mu$.
This is a generalization of the factorial powers.
We have $(a\u(k))=(a\u k)$ and $(a\u(1^k))=(a\f k)$.
Next put
$$
\nS_\mu= \frac1{(n\u\mu)}\S_\mu\quad\in\Zgln \,. \tag 2.15
$$
By virtue of \tht{2}, vanishing theorem, and the well known formula
for the dimension of the representation $\pi_\l$
$$
\dim_{GL(n)}\l = \frac{(n\u\mu)}{H(\mu)} \,, \tag 2.16
$$
we have
$$
\tr\pi_\l(\nS_\mu) = \dt_{\l\mu},\quad |\l|\le|\mu|\,. \tag 2.17
$$

\subhead 2.7 \endsubhead
There is a quite far analogy between central elements $\Sm$
and characters of a finite group $G$ considered as elements 
of the group algebra $\C[G]$ of $G$. The relationship
between quantum immanants and characters of the symmetric
groups is especially close (see sections 4 and 5 below).
For example, the vanishing and characterization theorems
should be considered as an analog of the
orthogonality relations for characters of groups.

Suppose $\pi$ and $\rho$ are two non-equivalent 
irreducible representations of $G$ and suppose $\chi^\pi$
and $\chi^\rho$ are the corresponding irreducible characters. The 
orhtogonality of $\chi^\pi$ and $\chi^\rho$
$$
\sum_{g\in G} \chi^\rho(g) \cdot \chi^\pi(g^{-1}) = 0
$$
can be rewritten as
$$
\pi\left(\,
\sum_{g\in G} \chi^\rho(g) \cdot g^{-1} \right) = 0\,,
$$
that is as vanishing of some central element in the 
representation $\pi$. In the algebra $\Ugln$ there are
no elements that vanish in all but one irreducible
representations; however the quantum immanants vanish
is as many representations as possible.
The vanishing and characterization theorems play the same
(and, perhaps, more imporatant) role in the combinatorics 
of $s^*$-functions as the orthogonality relations in the
combinatorics of $s$-functions (see \cite{OO}).

\head 
3. Proof of the main theorem
\endhead
In this section we shall prove the higher Capelli identities 
\tht{1.7}.
\subhead 3.1\endsubhead
Consider the following element
$$
R(u)=1+u\cdot(12) \in \k[S(2)][u] \,.
$$
It is called the $R$-matrix. Normally this $R$-matrix is denoted
by $\check R$, but we do not use other $R$-matrices. The
following equation can be verified by direct calculation:
$$
R(u-v) \cdot E(u)\otimes E(v) =
E(v)\otimes E(u) \cdot R(u-v) \,. \tag 3.1
$$
This is a version of the famous $RTT=TTR$ equation \cite{RTF}.
The equation \tht{1} is equivalent to the commutation
relations between the generators $\Eij$ of $\Ugln$.

\subhead 3.2\endsubhead
The second key fact we need is the Young orthogonal form \cite{JK}.
Put $s_i=(i,i+1)\in S(k)$. Consider the action of $s_i$ in
the Young orthogonal basis. Let $T$ be a standard tableau
and let $T'=s_i T$ (that is $T$ with $i$ and $i+1$
permuted).
Put
$$
r=c_T(i+1)-c_T(i)\,.
$$
If $T'$ is not a standard tableau then $r=\pm1$ and
$$
s_i \xi_T = \pm \xi_T\,.
$$
If $T'$ is a standard tableau then $|r|>1$ and
$$
s_i|_{\k\xi_T+\k\xi_{T'}} =
\left[
\matrix
r^{-1} & (1-r^{-2})^{1/2} \\
(1-r^{-2})^{1/2} & -r^{-1} 
\endmatrix
\right]\,. \tag 3.2
$$
Put $R_i(u)=1+u\cdot s_i$ and put
$$
R_i(T)=R_i(-r)\,.
$$
Clearly,
$$
R_i(T)|_{\k\xi_T+\k\xi_{T'}} =
\left[
\matrix
0 & -r(1-r^{-2})^{1/2} \\
-r(1-r^{-2})^{1/2} & 2 
\endmatrix
\right]\,, \tag 3.3
$$
if $T'$ is a standard tableau. In any case
$$
(R_i(T) \xi_T,\xi_T) = 0 \,. \tag 3.4
$$
Recall that $P_T$ is the orthogonal projection onto $\xi_T$.
It follows that if $T'$ is standard then
$$
R_i(T) P_T =
P_{T'} R_i(T) P_T  \tag 3.5
$$
and 
$$
P_{T'} = 
(r^2-1)^{-1} \,
R_i(T) P_T R_i(T) \,. \tag 3.6
$$
\subhead 3.3\endsubhead
Given a standard tableau $T$ put
$$
E(T)= E\otimes E(c_T(2))\otimes \dots \otimes E(c_T(k)) 
$$
Remark that \tht{1} reads as
$$
R_i(T)E(T)=E(T')R_i(T)\,, \tag 3.7
$$
where $T'=s_i T$. The third key fact we need is
\proclaim{\smc Proposition}
$$
E(T) P_T=P_T E(T) P_T \,. \tag 3.8
$$
\endproclaim
This proposition is apparently due to I.~V.~Cherednik \cite{Ch}.
(See also \cite{JKMO}).
\demo{Proof}
First suppose that $T$ is the row tableau $\Tr$ of shape $\mu$
that is tableau filled in from left to right from top
to bottom.
For example if $\mu=(3,2,1)$ then $\Tr$ looks as follows
$$
\Tr=
\matrix
1&2&3\\
4&5&\\
6
\endmatrix
$$
Denote by $\P$  the row symmetrizer of $\Tr$
$$
\P=\sum_{s \text{ preserves the rows of  } \Tr\,} s
$$
and by $\Q$ the column antisym\-me\-tri\-zer of $\Tr$
$$
\Q=\sum_{s \text{ preserves the columns of  } \Tr\,}
\sgn(s)\cdot s \,.
$$
The product 
$$
\P\Q \in \kSk \tag 3.9
$$
is known as the Young symmetrizer correspondinding
to the tableau $\Tr$.  
Denote by $W^\mu$ the irreducible $S(k)$-module labeled by $\mu$.
The Young symmetrizer acts as zero operator 
in all irreducible $S(k)$-modules except $W^\mu$.
It is clear that
$$
\P^2=\mu!\,\P\,,
$$
where $\mu!\,=\mu_1!\,\mu_2!\,\dots$. It is well known \cite{JK} that
$$
(\P \Q)^2
=H(\mu) \P \Q \,. \tag 3.10
$$
Hence the element
$$
\frac1{\mu!\,H(\mu)}\P \Q \P \in \kSk
$$
acts as an orthogonal projection in $W^\mu$ 
and as zero operator in other representations
of $S(k)$. Again it is well known that 
the vector $\xi_{\Tr}$ is the unique vector in $W^\mu$
that is invariant under the action of the row-stabilizer
of $\Tr$. Therefore
$$
P_{\Tr} = 
\frac1{\mu!\,H(\mu)}\P \Q \P \,. \tag 3.11
$$
By virtue of \tht{10} we have equality of right ideals
$$
P_{\Tr} \, \kSk = \P \Q \, \kSk \,. \tag 3.12
$$
Denote the right ideal \tht{12} by $I$. Consider the annihilator $J$
of $I$ in the semisimple algebra $\kSk$
$$
J=\{x\in\kSk\,|\, xI=0\}\,.
$$
This is a left ideal in $\kSk$.
Put
$$
M_j=\mu_1+\dots+\mu_j \,, \quad j=1,2,\dots \,.
$$
In \cite{JKMO} it is shown that $J$
is the left ideal generated by
the following elements $J_i$, $i=1,\dots,k-1$
$$
J_i=
\cases
1-(i,i+1), & i\ne M_1,M_2,\dots \\
(1+s_{M_{j-1}+1})
(1+2s_{M_{j-1}+2}) \dots
(1+\mu_j s_{M_j}), &i=M_j \,.
\endcases
$$
Now write
$$
E(\Tr)P_{\Tr}=P_{\Tr} E(\Tr) P_{\Tr} +
(1-P_{\Tr}) E(\Tr) P_{\Tr}\,.
$$
We are going to show that the second summand vanishes.
We have $(1-P_{\Tr})P_{\Tr}=0$ and hence 
$(1-P_{\Tr})\in J$. Therefore it suffices to check
that
$$
J_i E(\Tr) P_{\Tr}= 0, \quad i=1,2,\dots \,.
$$
Suppose $i\ne M_1,M_2,\dots $. Then
$$
J_i = R_i(\Tr) \,.
$$
Therefore by \tht{7}
$$
J_i E(\Tr) P_{\Tr} = E((\Tr)') J_i P_{\Tr} = 0 \,.
$$
If $i=M_j$ for some $j$ then we have to apply
the same relation \tht{7} several times. Hence \tht{8}
is proved for $\Tr$.

Next suppose for some $i$ both $T$ and $T'=s_i T$ 
are standard tableaux. Show that if \tht{8} holds for 
$T$ then it also holds for $T'$. Put
$$
r=c_T(i+1)-c_T(i)\,.
$$
If $T'$ is a standard tableau then $r\ne\pm1$. By \tht{6}
and \tht{7}
$$
\align
E(T')P_{T'} 
&= (r^2-1)^{-1} E(T') R_i(T) P_T R_i(T) \\
&= (r^2-1)^{-1} R_i(T) E(T) P_T R_i(T) \,.
\endalign
$$
By assumption this equals to
$$
(r^2-1)^{-1} R_i(T) P_T E(T) P_T R_i(T) \,.
$$
By \tht{5} this expression is stable under 
multiplication by $P_{T'}$ on the left. Thus
\tht{8} is proved for $T'$. This completes
the proof of the proposition. \qed
\enddemo

\subhead 3.4\endsubhead 
Recall that by definition
$$
\S_\mu=\tr E(T)P_T\,,
$$
where $T$ is a standard tableau of shape $\mu$.
\proclaim{\smc Proposition}
The quantum immanant $\S_\mu$ is a well defined
element of $\Ugln$. In other words, the trace
$$
\tr E(T)P_T
$$
does not depend on the choice of the standard
tableau $T$ of shape $\mu$.
\endproclaim
\demo{Proof}
Show that
$$
\tr E(T)P_T = \tr E(T')P_T' \,, \tag 3.14
$$
where $T'=s_i T$. Note that if $T'$ is a standard
tableau then $R_i(T)$ is invertible. 
Consider the following chain of equalities
$$
\alignat2
\tr E(T)P_T 
&= \tr R_i(T) E(T) P_T R_i(T)^{-1} &&\\
&= \tr E(T') R_i(T) P_T R_i(T)^{-1} \quad&&\text{by \tht{7}} \\
&= \tr E(T') P_{T'} R_i(T) P_T R_i(T)^{-1} \quad&&\text{by \tht{5}} \\
&= \tr P_{T'} E(T') P_{T'} R_i(T) P_T R_i(T)^{-1} \quad&&\text{by \tht{8}} \\
&= \tr E(T') P_{T'} R_i(T) P_T R_i(T)^{-1} P_{T'} \,. 
\endalignat
$$
It remains to prove that
$$
P_{T'} R_i(T) P_T R_i(T)^{-1} P_{T'} = P_{T'} \,. \tag 3.15
$$
Clearly \tht{15} holds up to a constant factor. Since
the highest terms of \tht{14} agree this constant equals 1.
The proposition is proved. \qed
\enddemo

\subhead 3.5\endsubhead
Recall that we identify the algebra of right-invariant differential operators
on $n\times n$ matrices with $\Ugln$. 
\proclaim{\smc Proposition} 
$$
\S_\mu,\Dm \in\Zgln\,. \tag 3.16
$$
\endproclaim
\demo{Proof}
Denote by $g_{ij}$ and $\gi_{ij}$ the matrix elements
of a matrix $g\in GL(n)$ and its inverse matrix $g^{-1}$.
The following equality is obvious
$$
\sum_k g_{ki} \gi_{jk} = \dt_{ij}\,. \tag 3.17
$$
Consider the adjoint action $\Ad(g)$ of $g$ in $\gln$
$$
\Ad(g)\cdot E_{ij} = \sum_{k,l} g_{ki} \gi_{jl} E_{kl} \,. \tag 3.18
$$
Under the adjoint action of $g$ the 
entries of the matrix $E(u)$ are
transformed as follows
$$
\alignat2
E(u) @>\Ad(g)>> 
& \sum_{i,j} 
\left(
\sum_{k,l} g_{ki} \gi_{jl} E_{kl} \right) \otimes
e_{ij} - u\sum_{i}1\otimes e_{ii} 
\quad&&\text{by \tht{18}} \\
&= \sum_{k,l} 
\left(
E_{kl} - u\dt_{kl} \right) \otimes
\left(\sum_{i,j} g_{ki} \gi_{jl} e_{ij} \right) 
\quad&&\text{by \tht{17}} \\
&= g' E(u) (g')^{-1} \tag 3.19
\endalignat
$$
The product \tht{19} is the product of the matrix $E(u)$
with entries in $\Ugln$ and two matrices with entries
in the ground field $\k$.
Consider the following element of $\Ugln$
$$
\tr(E(u_1)\otimes\dots\otimes E(u_k)\cdot s)\,,  \tag 3.20
$$
where $u_i\in\k$ and $s\in S(k)$ are arbitrary. By \tht{19} the adjoint action
of $g'$ takes \tht{20} to
$$
\tr(g^{\otimes k}
E(u_1)\otimes\dots\otimes E(u_k)
(g^{-1})^{\otimes k}
\cdot s) =
\tr(E(u_1)\otimes\dots\otimes E(u_k)\cdot s)\,.
$$
That is \tht{20} is an element of $\Zgln$.
In particular,
$$
\S_\mu \in \Zgln \,.
$$

Under the left action of an element $g\in GL(n)$
the matrices $X$ and $D$ are
transformed as follows
$$
X@>g>>g'X, \quad D@>g>>g^{-1}D \,.
$$
Therefore the left action of $g'$ takes $\Dm$ to
$$
\align
\Dm@>g>>
&\tr \big(g^{\otimes k} \cdot X^{\otimes k} \cdot (D')^{\otimes k} 
\cdot (g^{\otimes k})^{-1}
\cdot \chm/k! \big) \\
&=\tr \big(X^{\otimes k} \cdot (D')^{\otimes k}
\cdot \chm/k! \big) \\
&=\Dm \,.  
\endalign
$$
In the same way the $\Dm$ is invariant under the
right action of $GL(n)$. Therefore it represent an element
of $\Zgln$. \qed
\enddemo

\subhead 3.6\endsubhead
By definition put
$$
E(\mu)=\sum_T E(T)P_T \,, \tag 3.21
$$
where the summation is over all standard tableaux $T$ of shape $\mu$.
By proposition 3.4 we have
$$
\S_\mu = \frac 1{\dim\mu}\, \tr E(\mu) \,. \tag 3.22
$$
\proclaim{\smc Proposition}
$$
s E(\mu) = E(\mu) s, \quad \text{ for all } s\in S(k) \,.  \tag 3.23
$$
\endproclaim
\demo{Proof}
We can assume $s\in\kSk$. In $\kSk$ there is a basis of
matrix elements all irreducible representations of $S(k)$
corresponding to the Young orthonormal basis in each
representation. If
$s$ is a matrix element of a representation $\nu$, $\nu\ne\mu$,
then by \tht{8} both LHS and RHS of \tht{23} equal zero.
Suppose $s$ is a matrix element of the representation $\mu$.
The diagonal matrix elements in the Young basis are proportional
to $P_T$, where $T$ runs over standard tableaux of shape $\mu$.
Clearly
$$
P_T E(\mu) = E(\mu) P_T = E(T) P_T
$$
by \tht{8}. Suppose $s$ is a non-diagonal matrix element.
We can assume that $s$ takes $\xi_T$ to $\xi_{T'}$, $T'=s_i T$,
since such elements are generators. In this case $s$ is
proportional to 
$$
P_{T'} R_i(T) P_T \,.
$$
We have
$$
\alignat2
P_{T'} R_i(T) P_T E(\mu) 
&= P_{T'} R_i(T) P_T E(T) P_T &&\quad \text{by \tht{8}} \\
&= P_{T'} R_i(T) E(T) P_T &&\quad \text{by \tht{8}} \\
&= P_{T'} E(T') R_i(T) P_T &&\quad \text{by \tht{7}} \\
&= P_{T'} E(T') P_{T'} R_i(T) P_T &&\quad \text{by \tht{5}} \\
&= E(T') P_{T'} R_i(T) P_T &&\quad \text{by \tht{8}} \\
&= E(\mu) P_{T'} R_i(T) P_T \,,&&\quad  
\endalignat
$$
as desired. \qed
\enddemo

Suppose we have a sequence of indices 
$i_1\le i_2\le\dots\le i_k$ or 
$i_1\ge i_2\ge\dots\ge i_k$. Suppose that
exactly $\iota_1$ first indices are equal,
exactly $\iota_2$ next indices are equal and so on.
Then the stabilizer in $S(k)$ of the sequence
$(i_1,\dots,i_k)$ is isomorphic
to
$$
S(\iota)=S(\iota_1)\times S(\iota_2)\times\dots \,.
$$
We have
$$
|S(\iota)|=\iota!=\iota_1!\iota_2!\dots \,.
$$
For example, if all $i_j$ are distinct then $\iota=(1^k)$.
\proclaim{\smc Corollary}
$$
\S_\mu=\sum_{i_1\ge\dots\ge i_k} 1/{\iota!} 
\sum_T \sum_{s\in S(k)}
(s\cdot\xi_T,\xi_T) \, E_{i_1,i_{s(1)}} 
(E_{i_2,i_{s(2)}}-c_T(2)\dt_{i_2,i_{s(2)}}) \dots \,. \tag 3.24
$$
\endproclaim
\demo{Proof}
By \tht{23} all $k!/\iota!$ summands corresponding
to different rearrangements of $\{i_1,\dots,i_k\}$ make
the same contribution to the sum \tht{22}. \qed
\enddemo
In the same way we can write
$$
\S_\mu=\sum_{i_1\le\dots\le i_k} 1/{\iota!} 
\sum_T \sum_{s\in S(k)}
(s\cdot\xi_T,\xi_T) \, E_{i_1,i_{s(1)}} 
(E_{i_2,i_{s(2)}}-c_T(2)\dt_{i_2,i_{s(2)}}) \dots \,. \tag 3.24${}'$
$$
The formula \tht{24${}'$} turns into the original
definition \tht{1.1} of the Capelli element if $\mu=(1^n)$.

\subhead 3.7\endsubhead
Now we can calculate the eigenvalues of the 
quantum immanants.

By $\RTab(\mu,n)$ denote the
set of reverse column strict tableau $T$ of shape $\mu$ with
entries in $\{1,\dots,n\}$. By definition 
$T\in \RTab(\mu,n)$ if entries of $T$ weakly decrease along
the rows and strictly decrease along the columns. 
By definition, put
$$
\Sim(x_1,\dots,x_n)=
\sum_{T\in \RTab(\mu,n)\,} \, \prod_{\a\in\mu}(x_{T(\a)}-c(\a))\,, \tag 3.25
$$
where the product is over all cells $\a$ of $\mu$ and 
$c(\a)$ denotes the content of the cell $\a$. For factorial
Schur polynomials sums analogous  
to \tht{25} were considered by Biedenharn and Louck \cite{BL},
Chen and Louck \cite{CL}, Goulden and Hamel \cite{GH},
Macdonald \cite{M2} and others.

Since the content of the cell $(1,1)\in\mu$
equals $0$ the sum \tht{25} is stable in the 
following sense
$$
\Sim(x_1,\dots,x_n,0)=
\Sim(x_1,\dots,x_n)\,. \tag 3.26
$$
Let $\l$ be a partition.
By \tht{26} the number $\Sim(\l)$ is well defined.

\proclaim{\smc Proposition}
$$
\pi_\l(\S_\mu)=\Sim(\l)\,. \tag 3.27
$$
\endproclaim
\demo{Proof}
Apply \tht{24} to the highest vector. Since the highest
vector is annihilated by all $\Eij$ with $i<j$ we get
$$
\pi_\l(\S_\mu)=
\sum_{i_1\ge\dots\ge i_k} 1/{\iota!} 
\sum_T \sum_{s\in S(\iota)}
(s\cdot\xi_T,\xi_T) \, \l_{i_1} 
(\l_{i_2}-c_T(2)) \dots \,. \tag 3.28
$$
Here the summation in $s$ is over the stabilizer $S(\iota)\subset S(k)$
of $i_1,\dots,i_k$. 

Given a standard tableau $T$ denote by $\Tt=i(T)$ the tableau
obtained by replacing each number $j$ in $T$ by $i_j$.
The entries in rows and columns of $i(T)$ weakly decrease. We have
$$
\pi_\l(\S_\mu)=
\sum_{i_1\ge\dots\ge i_k} 1/{\iota!} 
\sum_{\Tt}
\l_{i_1} 
(\l_{i_2}-c_T(2)) \dots  
\left(
\sum_{s\in S(\iota)}
\sum_{T\in i^{-1}(\Tt)}
(s\cdot\xi_T,\xi_T) \right) \,. \tag 3.29
$$
Here $\Tt$ runs over all tableaux with entries $i_1,\dots,i_k$
and weakly decreasing rows and columns.
Next we show that
$$
\left(
\sum_{s\in S(\iota)}
\sum_{T\in i^{-1}(\Tt)}
(s\cdot\xi_T,\xi_T) \right) =
\cases 
\iota!, &\text{ if }\Tt\text{ is column strict}\\
0,&\text{ otherwise.} 
\endcases \tag 3.30
$$
Consider the diagram $\mu$ as disjoint union of
skew diagrams $\mu_1,\mu_2,\dots$ (which depend on the
tableau $\Tt$) defined as follows. The diagram
$\mu_1$ consists of first $\iota_1$ cells of $\Tt$,
the diagram $\mu_2$ consists of next $\iota_2$ cells
of $\Tt$ and so on. Then it is easy to see that
the LHS of \tht{30} equals 
$$
\prod_m \sum_{s\in S(\iota_m)}  \chi^{\mu_m}(s)\,. \tag 3.31
$$
It is an elementary fact from the representation
theory of the symmetric group that for any skew
diagram $\eta$
$$
\sum_{s\in S(|\eta|)}  \chi^{\eta}(s)=
\cases
|\eta|!, &\text{if }\eta\text{ is a horizontal strip,}\\
0,&\text{otherwise}\,. 
\endcases \tag 3.32
$$
This proves \tht{30} and therefore 
$$
\pi_\l(\S_\mu)=
\sum_{i_1\ge\dots\ge i_k} \,
\sum_{\Tt\in \RTab(\mu,n)\,} \,\, \prod_{\a\in\mu}(\l_{\Tt(\a)}-c(\a))\,.
$$
Finally observe that the summation over $i_1,\dots,i_k$
can be eliminated if we allow $\Tt$ to range over
all inverse column strict tableaux of shape $\mu$ and
entries $1,\dots,n$. Hence 
$$
\pi_\l(\S_\mu)=
\sum_{\Tt\in \RTab(\mu,n)\,} \,\, \prod_{\a\in\mu}(\l_{\Tt(\a)}-c(\a))\,. 
\tag 3.33
$$
Clearly this is the desired formula. \qed
\enddemo

\subhead 3.8\endsubhead
In this paragraph we prove the following 
\proclaim{\smc Proposition}
$$
\pi_\l(\S_\mu)=\sm(\l)\,. \tag 3.34
$$
\endproclaim 
By virtue of \tht{27} this proposition is equivalent
to the following {\it combinatorial} formula for
$s^*$-functions.
$$
\sm(x_1,x_2,\dots)=\Sim(x_1,x_2,\dots) \,. \tag 3.35
$$
This formula is eqivalent to the analogous formula
for factorial Schur functions (see \cite{CL},\cite{GG} and 
\cite{M2}). A proof of \tht{35} is
given also in \cite{OO}. Below we give one more
proof of \tht{35} based on the characterization
theorem.

\demo{Proof}
The number $\Sigma_\mu(\l)$ equals by \tht{27} to the
eigenvalue of an element of $\Zgln$ and hence 
$\Sigma_\mu\in\Ls(n)$.

Show that 
$\Sigma_\mu(\l)=0$ unless $\mu\in\l$.
Moreover, show that
$$
\prod_{\a\in\mu}(\l_{T(\a)}-c(\a)) = 0 \tag 3.36
$$
for all $T\in\RTab(\mu,n)$ unless $\mu\in\l$.
Denote the LHS of \tht{36} by $\Pi_T$. Put
$\l_{(i,j)}=\l_{T(i,j)}$.
Since $T\in\RTab(\mu,n)$ we have
$$
\l_{(1,1)}\le\l_{(1,2)}\le\dots\le\l_{(1,\mu_1)} \,. \tag 3.37
$$
If $\Pi_T(\l)\ne 0$ then
$$
\l_{(1,1)}\ne 0,\quad\l_{(1,2)}\ne1,\quad\l_{(1,3)}\ne2,\dots \tag 3.38
$$ 
By \tht{37} and \tht{38} we have
$$
\l_{(1,1)}\ge1,\quad\l_{(1,2)}\ge2,\quad\l_{(1,3)}\ge3,\dots \tag 3.39
$$ 
Again since $T\in\RTab(\mu,n)$ we have
$$
T(1,i)<T(2,i)<\dots<T(\mu'_i,i) \tag 3.40
$$
and we have also
$$
i\le\l_{(1,i)}\le\l_{(2,i)}\le\dots\le\l_{(\mu'_i,i)} \tag 3.41
$$
for all $i$. Observe that \tht{40} and \tht{41} yield $\l'_i\ge\mu'_i$.
Thus $\Pi_T(\l)\ne 0$ implies $\mu\subset\l$. By the
characterization theorem $\Sigma_\mu$ equals $\s_\mu$
up to a constant factor. In order to see that this
factor equals 1 we can either compare the highest
terms of $\s_\mu$ and $\Sigma_\mu$ or calculate 
explicitly the unique non-vanishing summand in $\Sigma_\mu$.
\qed
\enddemo

It is interesting to look at shifted analogs of
elementary and complete homogeneous functions. By \tht{35}
we have
$$
\align
s^*_{(1^k)}(x)&=\sum_{i_1<\dots<i_k} 
(x_{i_1}+k-1)\dots(x_{i_{k-1}}+1) x_{i_k}\\
s^*_{(k)}(x)&=\sum_{i_1\le\dots\le i_k} 
(x_{i_1}-k+1)\dots(x_{i_{k-1}}-1) x_{i_k}\,.
\endalign
$$

\subhead 3.9\endsubhead
Now we can complete the proof of the theorem. Consider the 
difference 
$$
\S_\mu -\Dm
\,. \tag 3.42 
$$
As explained in paragraph 1.4 this is an element 
of degree $<|\mu|$. Next \tht{42} is a central element.
Prove that it vanishes in all representations $\pi_\l$
such that $|\l|<|\mu|$. We have proved in the
previous paragraph that $\pi_\l(\S_\mu)=0$ for such $\l$. 
The differential operator vanishes also.
Indeed all irreducible $GL(n)$-submodules
of $\k[\M(n)]$ with highest weight $\l$
consist of
polynomials of degree $|\l|<|\mu|$.
Such polynomials are clearly
annihilated by the operator $\Dm$. Thus by 
the characterization theorem \tht{42}
equals zero. This concludes the proof of the theorem.

\subhead 3.10 \endsubhead
The quantum immanant $\S_\mu$ can be expressed in terms of 
the Young symmetrizer \tht{9}. We keep the
notations of paragraph 3.3. The element
$$
H(\mu)^{-1} \P \Q 
$$
is an idempotent proportional to the Young
symmetrizer. 
Consider the following element of $\Ugln$
$$
H(\mu)^{-1} \,\tr E(\Tr) \P \Q \,. 
\tag 3.43
$$
We have
$$
\alignat 2
H(\mu)^{-1} \,\tr E(\Tr) \P \Q 
&= H(\mu)^{-1} \, \tr E(\Tr) P_{\Tr} \P \Q  
&&\quad\text{by\tht{11}}\\
&= H(\mu)^{-1} \, \tr P_{\Tr} E(\Tr) \P \Q  
&&\quad\text{by\tht{8}}\\
&= H(\mu)^{-1} \, \tr E(\Tr) \P \Q P_{\Tr} &&\\
&= \tr E(\Tr) P_{\Tr}   
&&\quad\text{by\tht{11}} \\
&=\Sm\,. &&
\endalignat
$$

\head
4.Higher Capelli identities for Schur-Weyl duality.
\endhead

\subhead 4.1 \endsubhead
Consider the space of tensors $(\k^n)^{\otimes K}$.
It is a multiplicity free $GL(n)\times S(K)$-module;
so we can look for Capelli identities (in the sence
of \cite{HU}) for it. 

Suppose $k\le K$ and $|\mu|=k$. Embed $S(k)$ in $S(K)$.
Denote by $\Ind\chm$ the induced character of $S(K)$.
By the Frobenius formula
$$
\Ind\chm =\sum_{t\in S(K)/S(k)} t\cdot\chm\cdot t^{-1} \,, \tag 4.1
$$
in other words $\Ind\chm$ is proportional  to the averaging of
$\chm\in\kSk$ over the group $S(K)$.

Let $\tau$ denote the representation of the
group $GL(n)$ in the space 
$(\k^n)^{\otimes K}$ and let $\sigma$ denote the representation of
the group $S(K)$ in the same space.
\proclaim{\smc Theorem}
$$
\tau(\S_\mu)=\sigma(\Ind \chm\ov(K-k)!) \,.\tag 4.2
$$ 
\endproclaim
\demo{Proof}
Denote by $\M(n,K)$ the space of rectangular $n\times K$
matrices. Let $\{e_i\}$, $i=1,\dots,n$ be the 
standard basis of $\k^n$.
Embed $(\k^n)^{\otimes K}$ in $\k[\M(n,K)]$ as follows
$$
e_{i_1}\otimes\dots\otimes e_{i_K} \to
x_{i_1 1}\dots x_{i_K K} \,. \tag 4.3
$$
This embedding is $GL(n)$-equivariant. By \tht{1.7} the operator
$\tau(\S_\mu)$ becomes $\Dm$.
Consider the action of the group $S(K)$
$$
s\cdot x_{i_1 1}\dots x_{i_K K} 
= x_{i_1 s^{-1}(1)}\dots x_{i_K s^{-1}(K)} \,.
$$
By its very definition the operator $\Dm$ acts
as follows
$$
x_{i_1 1}\dots x_{i_K K} 
@>\Dm>>
\sum_{t\in S(K)/(S(k)\times S(K-k))}\quad
\sum_{s\in S(k)} \chm(s) \,
(tst^{-1}) \cdot
x_{i_1 1}\dots x_{i_K K} \,. 
$$
Thus $\Dm$ acts in the same way as $\Ind \chm\ov(K-k)!$
\qed
\enddemo

Another approach to Capelli-type identities for
Schur-Weyl duality was developed in \cite{KO}

\head
5. Further properties of quantum immanants.
\endhead

The results of this sections are from \cite{OO} (only
the proofs differ). This results are based on \tht{2.2};
that is they are essentially properties of $s^*$-functions.

The proofs below use higher Capelli identities. One
can take a short-cut and deduce all propositions 
directly from the  characterization theorem. Such
proofs can be found in \cite{OO}.

\subhead 5.1 \endsubhead
We have considered $\Sm\in\Zgln$ where $n$ was a
fixed number. Now let $n$ vary.

Suppose $N>n\ge\ell(\mu)$. Consider the composition of the two maps
$$
\Zgln \to \Cal U(\frak g\frak l(N)) \to 
\frak Z(\frak g\frak l(N))\,, \tag 5.1
$$
where the first arrow is the natural inclusion and the
second one is the $\frak{gl}(N)$-invariant projection.
If $\k=\C$ then this composition is the  averaging
over the group $U(N)$. We call this map
the {\it averaging} map. We denote the averaging
of $\xi\in\Zgln$ by $\la\xi\ra_N$.

In order to avoid confusion denote by 
by $\nS_{\mu|n}$
the normalized quantum $\mu$-immanant in $\Zgln$.
We call the following property the {\it coherence}
of quantum immanants.

\proclaim{\smc Proposition \cite{OO}}
$$
\la \nS_{\mu|n} \ra_N = 
\nS_{\mu|N} \tag 5.2
$$
\endproclaim
\demo{Proof}
Identify $\Ugln$ with distributions supported at $1\in\M(n)$.
By $s_\mu(D)$ denote the polynomial \tht{1.9} in variables
$\partial_{ij}$. 
\proclaim{\smc Lemma 1}
$$
(\Sm,\phi)=
\big[ 
s_\mu(D)\cdot\phi
\big] (1) \,. \tag 5.3
$$
\endproclaim
\demo{Proof of Lemma}
The higher Capelli identity \tht{1.7} asserts that
$$
\align
(\Sm,\phi)&=
\big[
\tr(X^{\otimes k}\cdot (D')^{\otimes k} \cdot
\chm\ov k!) \cdot \phi
\big] (1) \\
&=
\big[
\tr(
(D')^{\otimes k} \cdot
\chm\ov k!) \cdot \phi
\big] (1) \,, \tag 5.4
\endalign
$$
where we used the fact that $x_{ij}(1)=\dt_{ij}$.
By \tht{1.9} we have
$$
(\Sm,\phi)=
\big[ 
s_\mu(D')\cdot\phi
\big] (1) \,. 
$$
Since $\chm(s)=\chm(s^{-1})$ for all $s\in S(k)$ we
have 
$$
s_\mu(D')=s_\mu(D) \,. \qed
$$
\enddemo

As in \tht{1} consider the composition of the
following inclusion and projection
$$
S(\gln)^{GL(n)}\to
S(\glN)\to
S(\glN)^{GL(N)}\,, \tag 5.5
$$
where $S(\gln)$ is the symmetric algebra of $\gln$
and $S(\gln)^{GL(n)}$ denotes the invariants of the
adjoint action of $GL(n)$. We call this map the
averaging map also. Using the invariant scalar
product in $\gln$
$$
(A,B)=\tr AB, \quad A,B\in\gln\,, \tag 5.6
$$
we construct a similar averaging map
$$
\k[\M(n)]^{GL(n)}\to
\k[\M(N)]\to
\k[\M(N)]^{GL(N)}\,. \tag 5.7
$$
Again to avoid confusion denote $s_{\mu|n}(X)$
the polynomial \tht{1.9} in matrix elements of a
$n\times n$ matrix $X$.

\proclaim{\smc Lemma 2}
$$
\alignat2
\la s_{\mu|n}(X) \ra_N &= 
\frac{(n\u\mu)}{(N\u\mu)}
\,s_{\mu|N}(X)&\qquad &\in\k[\M(N)]^{GL(N)}\,, \tag 5.8\\
\la s_{\mu|n}(D) \ra_N &= 
\frac{(n\u\mu)}{(N\u\mu)}
\,s_{\mu|N}(D)&\qquad &\in S(\glN)^{GL(N)}\,. \tag 5.9
\endalignat
$$
\endproclaim
\demo{Proof of Lemma}
Recall that $s_\mu(X)$ equals the trace of $X=(x_{ij})$
in the irreducible $GL(n)$-module with highest
weight $\mu$. Consider the matrix element $f_\l$ 
corresponding to the highest vector 
$$
f_\l=
\prod_{i=1}^{\ell(\mu)}
\det\left[
\matrix
x_{11}&\hdots&x_{1i}\\
\vdots&&\vdots\\
x_{i1}&\hdots&x_{ii}
\endmatrix
\right]^{\mu_i-\mu_{i+1}}\,\in\k[\M(\ell(\mu))] \,.
$$
The averaging of a matrix element equals the trace:
$$
\la f_\l\ra_n = \frac1{\dim_{GL(n)}\l}\, s_{\mu|n}(X)\,, \tag 5.10
$$
where $\dim_{GL(n)}\l$ denotes the dimension of the
irreducible $GL(n)$-module with highest weight $\l$.
Hence
$$
\align
\la s_{\mu|n}(X) \ra_N 
&= \dim_{GL(n)}\l\,\la \, f_\l \ra_N \\
&=
\frac{\dim_{GL(n)}\l}{\dim_{GL(N)}\l}\,
s_{\mu|N}(X)\,. \tag 5.11
\endalign
$$
Now \tht{8} follows from \tht{2.16}.
The claim \tht{9} follows from \tht{8} and the formula \tht{6}
for the invariant scalar product.
\qed
\enddemo

The proposition follows immediately from \tht{2},\tht{9}
and the definition of $\nS_\mu$. \qed
\enddemo

\subhead 5.2 \endsubhead
There is a map in the inverse direction
$$
\ZglN\to\Zgln\,, 
$$
which is the restriction of invariant differtial
operators on $\M(N)$ 
to the invariant subspace $\k[\M(n)]$. This map
was studied by Olshanski. It plays the central role 
in \cite{Ol1}. It follows from the very
definition of the operator $\Dm$ that
$$
\S_{\mu|N}@>\text{ restriction }>>\S_{\mu|n} \tag 5.12
$$
provided $n\ge\ell(\mu)$. On the level
of eigenvalues \tht{12} is equivalent to the
stability \tht{2.4} of $s^*$-functions.

\subhead 5.3 \endsubhead
Suppose $|\l|=K$. By $\dim \l/\mu$ denote
the dimension of the skew Young diagram $\l/\mu$. This number 
equals
$$
\dim \l/\mu = \la\Res\chi^\l,\chm\ra_{S(k)}\,, \tag 5.14
$$
where $\Res\chi^\l$ is the restriction of $\chi^\l$ to
$S(k)$ and $\la\cdot,\cdot\ra$ is the standard scalar
product of functions on $S(k)$
$$
\la\phi,\psi\ra_{S(k)}=
(1/k!)\sum_{s\in S(k)} \phi(s) \psi(s^{-1}) \,.
$$
Equivalently $\dim\l/\mu$ is equal
to the number of paths from $\mu$ to $\l$ in the Young
graph. Recall that the Young graph is the oriented
graph whose vertices are partitions and two partitions
$\mu$ and $\nu$ are connected by an edge (write 
$\mu\nearrow\nu$) if $\nu/\mu$ is a single cell. 
\proclaim{\smc Proposition \cite{OO}}
$$
\frac{\dim \l/\mu}{\dim \l} =
\frac{s^*_\mu(\l)}{(|\l|\f |\mu|)} \,. \tag 5.15
$$  
\endproclaim
\demo{Proof}
Compare the eigenvalues of both sides of \tht{4.2}
in the irreducible submodule corresponding to $\l$.
By definition of $\S_\mu$ its eigenvalue equals $s^*_\mu(\l)$.
Calculate the eigenvalue of $\Ind\chm$. Its trace
equals
$$
\align
\chi^\l(\Ind\chm) 
&= K!\, \la\chi^\l,\Ind\chm\ra_{S(K)} \\
&= K!\, \la\Res\chi^\l,\chm\ra_{S(k)} 
\quad\text{by the Frobenius reciprocity }\\
&= K!\, \dim\l/\mu \,. \tag 5.16
\endalign
$$
Hence the eigenvalue of the left hand side of \tht{4.2} equals
$$
\frac{K!}{(K-k)!} 
\frac{\dim \l/\mu}{\dim \l} =
(|\l|\f |\mu) \frac{\dim \l/\mu}{\dim \l} \,.
$$
This yields \tht{15}. \qed
\enddemo

\subhead 5.4 \endsubhead
The main application of the formulas \tht{2} and \tht{15}
is the explicit solution of the
two following problems:
\roster
\item
given an element $s\in\kSk$ and a character $\chi^\l$ 
of the group $S(K)$, $K>k$ find
$$
\frac{\chi^\l(s)}{\dim\l}\,, 
$$
\item
given an element $\xi\in\Ugln$ and a representation 
$\pi_\l$ of the group $GL(N)$, $N>n$ find
$$
\frac{\tr \pi_\l(\xi)}{\dim\l} = \pi_\l(\la\xi\ra_N) \,,
$$
where $\la\xi\ra_N$ is the averaging of $\xi$.
\endroster

Indeed, it is clear that $s$ and $\xi$ can be assumed
to be central. In the center of $\kSk$ and $\Ugln$
there is the basis of irreducible characters and
quantum immanants respectively. Finally observe
that the problems are linear in $s$ and $\xi$
respectively.

These problems play the key role in the
ergodic method of A.~M.~Vershik and S.~V.~Kerov 
\cite{VK}. In fact the understanding of the
papers \cite{VK} was the original aim of 
G.~Olshanski and me. These problems are also
discussed in \cite{KO}. 

The solution of similar problems of other
classical groups will be given in a 
forthcoming paper by G.~Olshanski and me.

\enddocument

\end